\documentclass[executivepaper]{article}
\usepackage{amsmath,amssymb}
\usepackage{bigints}
\usepackage{indentfirst}
\usepackage{graphicx}
\usepackage{fancybox}
\usepackage{ascmac}
\usepackage{latexsym}
\usepackage{bm}
\usepackage{amsfonts}
\usepackage[all]{xy}
\usepackage{color}
\usepackage{tikz}
\usetikzlibrary{cd}
\usetikzlibrary{intersections, calc, arrows}
\usepackage{here}
\usepackage{amsthm}

\usepackage[implicit=false]{hyperref}

\makeatletter
\@addtoreset{equation}{section}
\makeatother

\makeatletter 
\newcommand\fleqnoff{\@fleqnfalse\@mathmargin\@centering} 
\newcommand\fleqnon[1][\leftmargini]{\@fleqntrue\@mathmargin=#1\relax 
\@ifundefined{mathindent}{\let\mathindent\@mathmargin}{}} 
\makeatother

\title{\textbf{\large THE DENSITY OF BORROMEAN PRIMES}}
\author{\normalsize Yuki ISHIDA, Atsuki KURAMOTO and Dingchuan ZHENG}
\date{}

\theoremstyle{definition}

\newtheorem{thm}{Theorem}[section]
\newtheorem*{thm*}{Theorem}
\newtheorem{dfn}[thm]{Definition}
\newtheorem{lem}[thm]{Lemma}
\newtheorem{prp}[thm]{Proposition}
\newtheorem{cor}[thm]{Corollary}
\newtheorem{ex}[thm]{Example}

\newcommand{\Thm}[3]{ 
\begin{thm}[#1]\label{#2}
 \textit{#3}
\end{thm}
\vspace{0.1in}}

\newcommand{\Thmm}[3]{ 
\begin{thm*}[#1]\label{#2}
 \textit{#3}
\end{thm*}
\vspace{0.1in}}

\newcommand{\Lem}[3]{
\begin{lem}[#1]\label{#2}
  \textit{#3}
 \end{lem}
 \vspace{0.1in}}

\newcommand{\Prp}[3]{
 \begin{prp}[#1]\label{#2}
  \textit{#3}
 \end{prp}
 \vspace{0.1in}}

\begin{document}
 
\maketitle

\subsection*{Abstract}
\footnote[0]{2020 Mathematics Subject Classification: 11R44, 11N45}
\footnote[0]{Key words: Arithmetic Topology, R{\'e}dei extension, Borromean Primes, Refined Chebotarev Density Theorem}
In this paper, we study an asymptotic distribution of sets of primes satisfying certain ``linking conditions'' in arithmetic topology, namely, conditions given by the Legendre and R{\'e}dei symbols among sets of primes. As our Main Theorem, we prove an asymptotic density formula for Borromean primes among all primes. For the proof, we use the effective Chebotarev density formula under the Generalized Riemann Hypothesis and explicit computations of discriminants of the number fields involved in R{\'e}dei's extension.

\subsection*{Introduction}
\textit{Borromean primes} are triples of primes, which are defined as arithmetic analogues of Borromean rings, following the analogies between primes and knots in arithmetic topology (\cite{morishita2011knots}). The study of such triples of primes goes back to the work of R{\'e}dei in 1939. In \cite{Redei1939}, R{\'e}dei attempted to generalize Gauss's genus theory and introduced a certain triple symbol [$p_1,p_2,p_3$] for primes $p_1$, $p_2$, $p_3$. This symbol may be regarded as a triple generalization of the Legendre symbol $(\frac{p_1}{p_2})$. The {\'e}dei symbol describes the decomposition law of $p_3$ in a certain dihedral extension over $\mathbb{Q}$ of degree 8, which is unramified outside $p_1$, $p_2$ (cf. \cite{amano2014redei} for the uniqueness of such an extension). After a long silence, Morishita (\cite{10.3792/pjaa.76.18}, \cite{ca278c2c793644c982af3cf61d29e18e}) interpreted the R{\'e}dei symbol as an arithmetic analogue of Milnor's triple linking number of a link (\cite{milnor1954link}) in his study of arithmetic topology. Note that the Legendre symbol may be interpreted as an arithmetic analogue of Gauss's linking number of two knots (\cite[Chapter 4]{morishita2011knots}). Borromean primes are then defined as triples of primes $p_1,p_2,p_3$ satisfying the following conditions:
\[
 p_i \equiv 1 \,\, {\rm mod} \, 4 \,\, (i = 1,2,3), \hspace{0.1in} (\frac{p_i}{p_j}) = 1 \,\, (1 \leq i \neq j \leq 3) \hspace{0.1in} {\rm and} \hspace{0.1in} [p_1,p_2,p_3] = -1.
\]
The purpose of this paper is to study the distribution of Borromean primes.\par
The study of asymptotic distribution of primes also goes back to Gauss, and it is viewed as an origin of the so called arithmetic statistics nowadays. Gauss predicted the following asymptotic formula
\[
\pi(x) \sim \frac{x}{\log x},
\]
where $f(x) \sim g(x)$ means that $f(x)/g(x) \rightarrow 1$ when $x \rightarrow +\infty$. This formula was proved independently by Hadamard and de la Vall{\'e}e Poussin. In 19th century, the following formula was shown: for coprime integers $m(>1)$ and $a$,
\begin{equation}\label{eqdir}
\pi(x;a,m) \sim \frac{1}{\varphi(m)} \cdot \frac{x}{\log x},
\end{equation}
where $\varphi(m)$ is the Euler function. It was generalized to the following more general density theorem, known as the Chebotarev density theorem: Let $M$ be a number field and $M'$ be a finite Galois extension of $M$. For $\sigma \in {\rm Gal}(M'/M)$ and any positive real number $x$, we define
\[
\begin{split}
 &\pi_{M'/M}(x;\sigma)\\
 &:=\#\left\{\mathfrak{p} \in S_M^0 \, \middle| \, \mathfrak{p}:{\rm unramified \,\, in}\,\, M', \, {\rm N}_M\mathfrak{p}<x, \, [\frac{M'/M}{\mathfrak{p}}]=C(\sigma)\right\},
 \end{split}
\]
where $S_M^0$, ${\rm N}_M\mathfrak{p}$, $[\frac{M'/M}{\mathfrak{p}}]$ and $C(\sigma)$ mean the set of prime ideals of $M$, the absolute norm of $\mathfrak{p} \in S_M^0$, the Artin symbol for $\mathfrak{p} \in S_M^0$ and the conjugacy class of $\sigma$ in ${\rm Gal}(M'/M)$, respectively. Then the Chebotarev density theorem asserts
\begin{equation}\label{eqche}
\pi_{M'/M}(x;\sigma) \sim \frac{\#C(\sigma)}{\#G} \cdot \frac{x}{\log x}.
\end{equation}\noindent
Note that Dirichlet's theorem (\ref{eqdir}) is interpreted as a special case of Chebotarev's theorem (\ref{eqche}) in the $m$-th cyclotomic field $\mathbb{Q}(\zeta_m)$.\par
In this paper, we investigate the problem of asymptotic distribution of sets of primes satisfying certain ``linking condition” in arithmetic topology, namely, conditions given by the Legendre symbols and the R{\'e}dei symbol among primes. This problem requires some refinements of the classical results mentioned above. Let $\pi_{{\rm Borr}}(x)$ denote the number of Borromean primes $\{p_1,p_2,p_3\}$ with $p_i<x$ for $i=1,2,3$.  We prove the following asymptotic formula under the Generalized Riemann Hypothesis (abbreviated as GRH below).

\Thmm{Theorem \ref{mainthm} below}{}{If GRH is true, we have
\[
\lim_{x \to +\infty}\frac{\pi_{{\rm Borr}}(x)}{\#\{\{p_1,p_2,p_3\} \, | \, p_i < x \, ({\rm  for} \, \, i=1,2,3), \, p_i \neq p_j \, (i \neq j)\}}=\frac{1}{128}.
\]}
\noindent
We note that our theorem is not obtained by a straightforward application of the above classical density theorems and we need more elaborate analyses on the error term of the Chebotarev density theorem. For the proof, we employ the effective version of the Chebotarev density theorem under GRH (\cite{GRENIE2019441}) and, in order to apply the formula to our situation,  we compute explicitly the discriminants of the fields involved in R{\'e}dei's  extension. \par
This paper is organized as follows. In Section 1, we recall R{\'e}dei's extension, the R{\'e}dei triple symbols and Borromean primes. In addition, we evaluate the absolute discriminants of some extensions of R{\'e}dei's extension or its subextension. This value is necessary for the proof of our Main Theorem in Section 3. In Section 2, we calculate the density of pairs of distinct primes $\{p_1,p_2\}$ satisfying $p_1, p_2 \equiv 1$ mod 4, $(\frac{p_1}{p_2})=1$ and $p_1,p_2<x$ among all pairs of distinct primes $p_1,p_2<x$ (Theorem \ref{residenthm}). This value is also contributory to the proof of our theorem. In Section 3, we present our Main Theorem (Theorem \ref{mainthm}) and prove it.

\subsection*{Notations}
For integers $x,a$, and $b$, $x \equiv a \, (b)$ means that $x$ is congruent to $a$ modulo $b$. The symbol $(\frac{q}{p})$ and $[\frac{L/K}{\mathfrak{p}}]$ stands for the Legendre symbol and the Artin symbol, respectively. \par
Let $M$ be a number field and $M'$ be a finite extension of $M$. We write as $n_M$ the degree of the extension of $M$ over $\mathbb{Q}$. We denote by $\Delta_M$ the absolute discriminant of $M$ and by $\Delta_{M'/M}$ the relative discriminant of $M'$ over $M$. \par
For any positive real number $x$, $\pi(x)$ stands for the number of primes less than $x$. For such $x$ and any integers $a, \, b$, we denote by $\pi(x;a,b)$ the number of primes less than $x$ and congruent to $a$ modulo $b$.\par

\subsection*{Acknowledgements}
We would like to express our sincere gratitude to our advisor,  Professor Masanori Morishita, for his invaluable guidance and support throughout this study. We also extend our thanks to Professor Toshiki Matsusaka for sharing an analytical method, provided by Yuta Suzuki, for computing the density which is contained in Section 2. We also thank Deng Yuqi for useful discussions.

\section{R{\'e}dei's Dihedral Extension, R{\'e}dei Symbol and Borromean Primes}\label{fieldp1p2}

In this section, we recall the R{\'e}dei symbol defined in \cite{Redei1939}. We construct an extension of $\mathbb{Q}$ known as R{\' e}dei's extension. This is a Galois extension of degree 8 over $\mathbb{Q}$, where the Galois group is isomorphic to the dihedral group $D_4$ with order 8 and only specific primes ramify. Then we define the R\'{e}dei Symbol and introduce the Borromean primes. In addition, we calculate  the absolute discriminants of the fields obtained by adjoining $\sqrt{-1}$ to R{\'e}dei's extension or its subextension. These values are necessary later in the proof of our Main Theorem, in Section 3.\par
Let $p_1,p_2$ be distinct primes satisfying the conditions $p_1\equiv p_2\equiv 1 \,\, (4)$ and $(\frac{p_2}{p_1})=1$. In this case, the equation for $X,Y,Z$
\begin{equation}\label{redeieq}
X^2-p_1Y^2-p_2Z^2=0
\end{equation}
has non-trivial integer solutions $(x,y,z)$ (i.e. $xyz\neq 0$) satisfying the conditions
\begin{equation}\label{solcon}
 y \equiv 0 \, (2), \hspace{0.2in} x-y \equiv 1 \, (4)
\end{equation}
(\cite{Redei1939}). We choose one of such solutions and write it as $(x_0, y_0, z_0)$. We define
\[
 \alpha_2=x_0+y_0\sqrt{p_1}, \hspace{0.2in} \bar{\alpha}_2=x_0-y_0\sqrt{p_1},
\]
\[
 \alpha_1=\alpha_2+\bar{\alpha}_2+2z_0\sqrt{p_2}, \hspace{0.2in} \bar{\alpha}_1=\alpha_2+\bar{\alpha}_2-2z_0\sqrt{p_2},
\]
\[
 k_1=\mathbb{Q}(\sqrt{p_1}), \hspace{0.2in} k_2=\mathbb{Q}(\sqrt{p_2})
\]
and $k=\mathbb{Q}(\sqrt{p_1},\sqrt{p_2},\sqrt{\alpha_2})$. We consider the following tower of field extensions:
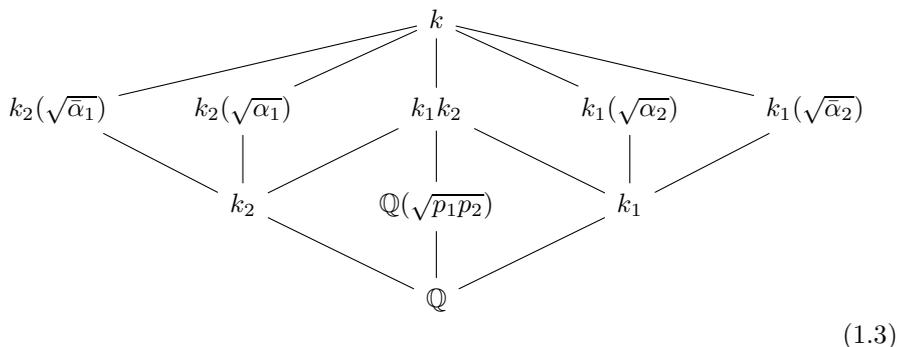
\begin{equation}\label{extredei}
\begin{tikzcd}
   & & k \\
 k_2(\sqrt{\bar{\alpha}_1}) \arrow[rru, -] & k_2(\sqrt{\alpha_1}) \arrow[ru, dash] & k_1k_2 \arrow[u, dash] & k_1(\sqrt{\alpha_2}) \arrow[lu, dash] & k_1(\sqrt{\bar{\alpha}_2}) \arrow[llu, dash] \\
 & k_2 \arrow[lu, dash] \arrow[u, dash] \arrow[ru, dash] & \mathbb{Q}(\sqrt{p_1p_2}) \arrow[u, dash] & k_1 \arrow[lu, dash] \arrow[u, dash] \arrow[ru, dash] \\
 && \mathbb{Q} \arrow[lu, dash] \arrow[u, dash] \arrow[ru, dash]
\end{tikzcd}
\end{equation}
It can be seen that the extension $k/\mathbb{Q}$ in the tower (\ref{extredei}) is a Galois extension whose Galois group is isomorphic to the dihedral group $D_4$. The field $k$ has the following uniqueness. 

\Prp{{\cite[Theorem 2.1, p3]{amano2014redei}}}{}{The field $k$ can be defined independently from the choice of the solution $(x,y,z)$ of the equation (\ref{redeieq}) satisfying conditions (\ref{solcon}) and it is characterized by the following two properties:
\begin{enumerate}
\item $k/\mathbb{Q}$ is a Galois extension with the Galois group being isomorphic to the dihedral group $D_4$. 
\item Primes ramified in $k/\mathbb{Q}$ are only $p_1$ and $p_2$ and both these ramification indices in $k/\mathbb{Q}$ are equal to $2$.
\end{enumerate}}
\noindent
For a prime $p_3$ different from $p_1$ and $p_2$ which satisfy the conditions
\[
 p_3 \equiv 1 \, (4) \hspace{0.2in} {\rm and} \hspace{0.2in} (\frac{p_1}{p_3})=(\frac{p_2}{p_3})=1,
\]
we define the symbol $[p_1,p_2,p_3]$ as follows:
\[
 [p_1,p_2,p_3]=\left\{\begin{array}{ll}
 1 & (p_3 \,\, {\rm is} \,\, {\rm completely} \,\, {\rm decomposed} \,\, {\rm in} \,\, k), \\
 \\
 -1 & (p_3 \,\, {\rm is} \,\, {\rm not} \,\, {\rm completely} \,\, {\rm decomposed} \,\, {\rm in} \,\, k).
 \end{array}\right.
\]
This definition is known to be well-defined regardless of the choice of $\alpha_2$ (\cite{Redei1939}). The symbol $[p_1, p_2, p_3]$ is called \textit{R{\'e}dei symbol}.\par
The Legendre symbol and the R{\'e}dei symbol can be regarded as arithmetic analogues of the linking number and Milnor's triple linking number, respectively,  in view of the analogies between knots and primes (\cite[Chapter 4, Chapter 8]{morishita2011knots}). Therefore following after the Borromean rings, \textit{Borromean primes} are defined as triples of primes $p_1,p_2,p_3$ satisfying the following conditions:
\[
 p_i \equiv 1 \,\, (4) \,\, (i = 1,2,3), \hspace{0.1in} (\frac{p_i}{p_j}) = 1 \,\, (1 \leq i \neq j \leq 3) \hspace{0.1in} {\rm and} \hspace{0.1in} [p_1,p_2,p_3] = -1.
\] 
The R{\'e}dei symbol has the following reciprocity. 

\Prp{{\cite[(37), p21]{Redei1939}}, {\cite[Theorem 3.2, p4]{amano2014redei}}}{reciredei}{For any permutation $\sigma \in \mathfrak{S}_3$ (the symmetric group of degree 3), we have the following reciprocity property:
\[
[p_{\sigma(1)},p_{\sigma(2)},p_{\sigma(3)}]=[p_1,p_2,p_3].
\]}
\noindent
Because of this proposition, it can be seen that the definition of Borromean primes is independent from the order of primes $p_1,p_2,p_3$. In other words, the definition is only due to the set of primes \{$p_1,p_2,p_3$\} satisfying the above condition.\par
Furthermore, we consider the following sequence of field extensions.
\begin{equation}\label{extp1p2}
\begin{tikzcd}
   & k(\sqrt{-1}) \\
 k \arrow[ru, dash]& (k_1k_2)(\sqrt{-1}) \arrow[u, dash]\\
 k_1k_2 \arrow[ru, dash] \arrow[u, dash] & \mathbb{Q}(\sqrt{-1}) \arrow[u, dash]\\
 \mathbb{Q} \arrow[ru, dash] \arrow[u, dash]
\end{tikzcd}
\end{equation}
In this case, all the extensions that appear in (\ref{extp1p2}) are finite Galois extensions. The following degrees of extensions are required for later computations in Section 3:
\[
 [k_1k_2:\mathbb{Q}]=4, \hspace{0.1in} [k:\mathbb{Q}]=8, \hspace{0.1in} [(k_1k_2)(\sqrt{-1}):\mathbb{Q}]=8, \hspace{0.1in} [k(\sqrt{-1}):\mathbb{Q}]=16.
\]\par
The absolute discriminants of $(k_1k_2)(\sqrt{-1})$ and $k(\sqrt{-1})$ will be used later to prove our Main Theorem. We state them here as the next propsitioin.

\Prp{}{discrilemma}{We have
\[
 \Delta_{(k_1k_2)(\sqrt{-1})}=2^8p_1^4p_2^4 \hspace{0.2in} {\it and} \hspace{0.2in} \Delta_{k(\sqrt{-1})}=2^{16}p_1^8p_2^{8}.
\]}

\begin{proof}
Let $K$ and $L$ be number fields with ${\rm GCD}(\Delta_K, \Delta_L)=1$, $[K : \mathbb{Q}]=m$ and $[L : \mathbb{Q}]=n$. Then it follows that
\begin{equation}\label{dispro}
\Delta_{KL} = \Delta_K^n \Delta_L^m
\end{equation}
(\cite[Theorem 5.13, p215]{mollin1999algebraic}). Note that direct computation leads to 
\[
\Delta_{k_1}=p_1, \hspace{0.1in} \Delta_{k_2}=p_2 \hspace{0.1in} {\rm and} \hspace{0.1in} \Delta_{\mathbb{Q}(\sqrt{-1})} = -2^2.
\]
Substituting $k_1$ and $k_2$ to $K$ and $L$ in the formula (\ref{dispro}) respectively, we have $\Delta_{k_1k_2}=p_1^2p_2^2$. Hence we obtain 
\[
 \Delta_{k_1k_2(\sqrt{-1})} = (p_1^2p_2^2)^2 \cdot (-2^2)^4 =2^8p_1^4p_2^4
\]
by substituting $k_1k_2$ and $\mathbb{Q}(\sqrt{-1})$ to $K$ and $L$ in (\ref{dispro}) respectively. \par
Next, we compute $\Delta_{k(\sqrt{-1})}$. In order to do so, we use the formula (\ref{dispro}) again. First, we evaluate $\Delta_k$ by the following formula: When $M$ is an intermediate field of $L/K$, it follows that
\begin{equation}\label{disnorm}
 \Delta_{L/K} = {\rm N}_{M/K}(\Delta_{L/M}) \cdot \Delta_{M/K}^{[L:M]}
\end{equation}
(\cite[Lemma 5.5, p198]{mollin1999algebraic}) where ${\rm N}_{M/K}(\Delta_{L/M})$ stands for the relative norm of $\Delta_{L/M}$ over $K$. Substituting $\mathbb{Q}, k_1$ and $k$ to $K,M$ and $L$ respectively, we have
\begin{equation}\label{eq13}
\Delta_{k/\mathbb{Q}}={\rm N}_{k_1/\mathbb{Q}}(\Delta_{k/k_1}) \cdot \Delta_{k_1/\mathbb{Q}}^4.
\end{equation}
Since $\Delta_{k_1/\mathbb{Q}}$ is equal to the ideal $(p_1)$, computing $\Delta_{k/k_1}$ is sufficient to obtain $\Delta_{k/\mathbb{Q}}$. Here is a proposition with relative discriminants proved by R{\'e}dei.

\Lem{{\cite[p6]{Redei1939}}}{redeilem}{The following holds:
\begin{enumerate}
\item[(1)] $\Delta_{k/k_1}=\Delta_{k_1k_2/k_1}{\rm N}_{k_1/\mathbb{Q}}(\Delta_{k_1(\sqrt{\alpha_2})/k_1})$.
\item[(2)] ${\rm N}_{k_1/\mathbb{Q}}(\Delta_{k_1(\sqrt{\alpha_2})/k_1}) = (p_2)$.
\end{enumerate}}
\noindent
Because of this Lemma, we only have to compute $\Delta_{k_1k_2/k_1}$. Since it is easily shown that the only prime $p_2$ is ramified in $k_1k_2/k_1$, $\Delta_{k_1k_2/k_1}$ is equal to some power of ideal $(p_2)$ in $k_1$. Using the formula (\ref{disnorm}) with $K=\mathbb{Q}$, $M=k_1$ and $L=k_1k_2$, we have $\Delta_{k_1k_2/k_1}=(p_2)$. According to Lemma \ref{redeilem}, we have
\[
\Delta_{k/k_1}=(p_2) \cdot (p_2) = (p_2)^2.
\]
Substituting this and $\Delta_{k_1/\mathbb{Q}}=(p_1)$ to (\ref{eq13}), we obtain
\[
\Delta_{k/\mathbb{Q}}={\rm N}_{k/\mathbb{Q}}((p_2)^2) \cdot (p_1)^4=(p_1^4p_2^4),
\]
so the absolute discriminant $\Delta_k$ is equal to $p_1^4p_2^4$. Since $p_1 \equiv p_2 \equiv 1 \,\, (4)$, $\Delta_k$ and $\Delta_{\mathbb{Q}(\sqrt{-1})}$ are relatively prime. Hence we have the following by the formula (\ref{dispro}) with $K=k$ and $L=\mathbb{Q}(\sqrt{-1})$:
\[
\Delta_{k(\sqrt{-1})}=(p_1^4p_2^4)^2 \cdot (-2^2)^8=2^{16}p_1^8p_2^8.
\]
This completes the proof of Proposition \ref{discrilemma}.
\end{proof}

\section{The Density of Quadratic Residue Primes}
In this section, we evaluate the density of pairs of primes each of which is congruent to $1$ modulo $4$ and quadratic residue modulo the other. This fact can be proved by using Theorem \ref{EChebo} in Section 3, but here we show it in another unconditional way. This way is due to Y. Suzuki.\par
In the following arguments, $p_i$ means a prime for all $i$. 

\Thm{}{residenthm}{We have the following asymptotic formula:
\[
\lim_{x \to +\infty}\frac{\#\{(p_1,p_2) \, | \, p_1,p_2<x, \, p_1 \equiv p_2 \equiv 1 \,(4), \, (\frac{p_1}{p_2})=1\}}{\pi(x)^2}=\frac{1}{8}.
\]}

\begin{proof}[\textbf{Proof} {\rm (Y. Suzuki)}]
We have
\begin{equation}\label{eq31}
\begin{split}
\#\{ &(p_1,p_2) \, | \, p_1,p_2 < x, p_1 \equiv p_2 \equiv 1(4), (\frac{p_1}{p_2}) = 1 \} \\
&=\frac{1}{2}\sum_{\substack{p_1,p_2 < x\\p_1,p_2 \equiv 1(4)}}\left(1+(\frac{p_1}{p_2})\right)-\frac{1}{2}\sum_{\substack{p < x\\p \equiv 1(4)}}1\\
&=\frac{1}{2}\pi (x ; 1, 4)^2+\frac{1}{2}\sum_{\substack{p_1,p_2 < x\\p_1,p_2 \equiv 1(4)}}(\frac{p_1}{p_2})-\frac{1}{2} \pi (x ; 1, 4).
\end{split}
\end{equation}
We evaluate the second term in the last RHS. We define
\[\begin{split}
a_m &=
\begin{cases}
1 & (\text{if $m$ : prime, $m \equiv 1(4)$}), \\
0 & (\text{otherwise})
\end{cases}\\
{\rm and} \hspace{0.2in} b_m &=
\begin{cases}
\displaystyle
\sum_{\substack{p_1 < x\\p_1 \equiv 1(4)}}(\frac{p_1}{m}) & \text{(if $m$ : square-free, odd),}\\
0 & \text{(otherwise),}
\end{cases}
\end{split}
\]
where $(\frac{n}{m})$ stands for the Jacobi symbol in the definition of $b_m$. Then we can write the second term in the last RHS of (\ref{eq31}) as
\[
E(x):=\sum_{\substack{p_2 < x\\p_2 \equiv 1(4)}}\sum_{\substack{p_1 < x\\p_1 \equiv 1(4)}}(\frac{p_1}{p_2})=\sum_{m < x}a_mb_m.
\]
By the Cauchy-Schwarz's inequality, we get
\begin{equation}\label{eq32}
\begin{split}
E(x)^2 = \left(\sum_{m < x}a_mb_m\right)^2 &\leq \left (\sum_{m < x}a_m^2\right) \left (\sum_{m < x}b_m^2\right)\\
&= \pi (x ; 1, 4) {\sum_{m < x}}^*\left(\sum_{\substack{p_1 < x\\p_1 \equiv 1(4)}}(\frac{p_1}{m})\right)^2\\
&= \pi (x ; 1, 4) {\sum_{m < x}}^*\left|{\sum_{n < x}}^*a_n \cdot (\frac{n}{m})\right|^2.
\end{split}
\end{equation}
where $\sum^*$ indicates summation over positive odd square-free integers only. We now employ the following theorem.

\Thm{{\cite[Theorem 1]{D1995}}, {\cite[Theorem 7.20]{MR2061214}}}{}{Let $M, N$ be positive integers, and let $a_1, \dots , a_n$ be arbitrary complex numbers. Then
\[
{\sum_{m \le M}}^*\left|{\sum_{n \le N}}^*a_n(\frac{n}{m})\right|^2 \le C_{\varepsilon} (MN)^{\varepsilon}(M+N){\sum_{n \le N}}^*\left|a_n\right|^2
\]
holds for any $\varepsilon > 0$, where $C_{\varepsilon}$ is constant depends only on $\varepsilon$.}

According to this theorem with $M=N=\lfloor x \rfloor$, where $\lfloor x \rfloor$ is the largest integer not exceeding $x$, we have for any real number $\varepsilon>0$
\[\begin{split}
E(x)^2 &\le C_{\varepsilon} \pi (x ; 1, 4)\lfloor x \rfloor^{1+2\varepsilon}{\sum_{n \le x}}^*|a_n|^2\\
&\leq C_{\varepsilon} \pi (x ; 1,4)^2x^{1+2\varepsilon}.
\end{split}
\]
Taking the square root of the both sides of this inequality, we obtain
\[
E(x) \leq \sqrt{C_{\varepsilon}} \cdot \pi(x;1,4)x^{\frac{1}{2}+\varepsilon}.
\]
Because of $\pi (x ; 1,4) \sim \frac{x}{2\log{x}}$ and Prime Number Theorem, we see that
\[
 \lim_{x \to +\infty}\frac{E(x)}{\pi(x)^2}=0.
\]
when we take $\varepsilon$ less than $\frac{1}{2}$. Therefore from the formula (\ref{eq31}), the following holds:
\[
\lim_{x \to +\infty}\frac{\#\{ (p_1,p_2) \, | \, p_1,p_2 < x, p_1 \equiv p_2 \equiv 1(4), (\frac{p_2}{p_1}) = 1 \}}{\pi (x)^2} = \frac{1}{8}.
\]
This is the desired result for Theorem \ref{residenthm}.
\end{proof}

\section{Main Theorem and Proof}

In the following argument, $p, \, p_i$ are primes for all $i$. We henceforth abbreviate the condition $p_i<x$ for $i=1,2,3$ as $p_i<x$ and similarly for other conditions with $p_i$. \par
Let $M$ be a number field and $M'$ be a finite Galois extension of $M$. Recall that for $\sigma \in {\rm Gal}(M'/M)$ and any positive real number $x$, we define
\[
\begin{split}
 &\pi_{M'/M}(x;\sigma)\\
 &:=\#\left\{\mathfrak{p} \in S_M^0 \, \middle| \, \mathfrak{p}:{\rm unramified \,\, in}\,\, M', \, {\rm N}_M\mathfrak{p}<x, \, [\frac{M'/M}{\mathfrak{p}}]=C(\sigma)\right\},
 \end{split}
\]
where $S_M^0$, ${\rm N}_M\mathfrak{p}$, $[\frac{M'/M}{\mathfrak{p}}]$ and $C(\sigma)$ mean the set of prime ideals of $M$, the absolute norm of $\mathfrak{p} \in S_M^0$, the Artin symbol for $\mathfrak{p} \in S_M^0$ and the conjugacy class of $\sigma$ in ${\rm Gal}(M'/M)$, respectively. \par
In this section, we prove our Main Theorem under GRH.
\Thm{Main Theorem}{mainthm}{If GRH is true, we have
\begin{equation}\label{maineq}
 \lim_{x \to +\infty}\frac{\pi_{{\rm Borr}}(x)}{\#\{\{p_1,p_2,p_3\} \, | \, p_i < x, \, p_i \neq p_j\}}=\frac{1}{128}.
\end{equation}}
\noindent
Recall that triples of primes $\{p_1,p_2,p_3\}$ are called Borromean primes when $p_i < x$, $p_i \equiv 1 \,(4)$, $(\frac{p_j}{p_i})=1$ and $[p_1,p_2,p_3]=-1$ and $\pi_{{\rm Borr}}(x)$ stands for the number of sets of Borromean primes $\{p_1,p_2,p_3\}$ with $p_i<x$.\par
We will prove Main Theorem by the following 3 steps:

\begin{enumerate}
\item The LHS of (\ref{maineq}) is transformed to (\ref{decomeq}) below in order to use the refined Chebotarev's density theorem (Theorem \ref{EChebo} below).
\item Some limits are derived by the refined Chebotarev (Proposition \ref{plimthm} below). 
\item Main Theorem is proved by using the proposition derived in the second step.
\end{enumerate}

\vspace{0.05in}

\begin{proof}[\textbf{Proof of Theorem \ref{mainthm}}]\quad\par
(step 1): Now, we calculate the LHS of (\ref{maineq}). The conditions for $p_1,p_2,p_3$
\[
 p_i < x, \hspace{0.2in} p_i \neq p_j
\]
are symmetric with $p_1,p_2,p_3$, that is, the same holds, even if indices are permuted. Then we have
\[
 \#\{\{p_1,p_2,p_3\} \, | \, p_i < x, \, p_i \neq p_j\}=\frac{1}{3!}\#\{(p_1,p_2,p_3) \, | \, p_i < x, \, p_i \neq p_j\}
\]
where $(*,*,*)$ denotes a vector with three components. Since the condition $[p_1,p_2,p_3]=-1$ is symmetric by Proposition \ref{reciredei}, the conditions
\[
 p_i < x, \, p_i \equiv 1 \,(4), \, (\frac{p_j}{p_i})=1, \, [p_1,p_2,p_3]=-1
\]
are also symmetric. Then
\[
 \begin{split}
 \pi_{{\rm Borr}}(x) &=\#\{\{p_1,p_2,p_3\} \, | \, p_i < x, \, p_i \equiv 1 \,(4), \, (\frac{p_j}{p_i})=1, \, [p_1,p_2,p_3]=-1\} \\
 &= \frac{1}{3!}\#\{(p_1,p_2,p_3) \, | \, p_i < x, \, p_i \equiv 1 \,(4), \, (\frac{p_j}{p_i})=1, \, [p_1,p_2,p_3]=-1\}
 \end{split}
\]
holds. Hence we derive
\begin{equation}\label{eqqqqq}
 \begin{split}
 &\frac{\pi_{{\rm Borr}}(x)}{\#\{\{p_1,p_2,p_3\} \, | \, p_i < x, \, p_i \neq p_j\}} \\
 &= \frac{\#\{(p_1,p_2,p_3) \, | \, p_i < x, \, p_i \equiv 1 \,(4), \, (\frac{p_j}{p_i})=1, \, [p_1,p_2,p_3]=-1\}}{\#\{(p_1,p_2,p_3) \, | \, p_i < x, \, p_i \neq p_j\}}.
 \end{split}
\end{equation}
Since direct computation yields
\[
 \begin{split}
 \pi(x)^3&=\#\{(p_1,p_2,p_3) \, | \, p_i<x\} \\
 &=\#\{(p_1,p_2,p_3) \, | \, p_i<x, \, p_i \neq p_j\}+3\pi(x)^2-2\pi(x),
 \end{split}
\]
we have
\[
 \begin{split}
 &\lim_{x \to +\infty}\frac{\#\{(p_1,p_2,p_3) \, | \, p_i<x, \, p_i \neq p_j\}}{\pi(x)^3} = \lim_{x \to +\infty}\left(1-\frac{3}{\pi(x)}-\frac{2}{\pi(x)^2} \right) = 1.
 \end{split}
\]
By this identity, the following holds:
\[
 \begin{split}
 &\lim_{x \to +\infty}\frac{\#\{(p_1,p_2,p_3) \, | \, p_i < x, \, p_i \equiv 1 \,(4), \, (\frac{p_j}{p_i})=1, \, [p_1,p_2,p_3]=-1\}}{\#\{(p_1,p_2,p_3) \, | \, p_i < x, \, p_i \neq p_j\}} \\
 &=\lim_{x \to +\infty}\frac{\#\{(p_1,p_2,p_3) \, | \, p_i < x, \, p_i \equiv 1 \,(4), \, (\frac{p_j}{p_i})=1, \, [p_1,p_2,p_3]=-1\}}{\pi(x)^3} \\
 &\hspace{2.5in} \times \lim_{x \to +\infty} \frac{\pi(x)^3}{\#\{(p_1,p_2,p_3) \, | \, p_i < x, \, p_i \neq p_j\}} \\
 &=\lim_{x \to +\infty}\frac{\#\{(p_1,p_2,p_3) \, | \, p_i < x, \, p_i \equiv 1 \,(4), \, (\frac{p_j}{p_i})=1, \, [p_1,p_2,p_3]=-1\}}{\pi(x)^3}.
 \end{split}
\]
According to this and (\ref{eqqqqq}), it is sufficient to prove Main Theorem that we evaluate the value
\[
 \lim_{x \to +\infty}\frac{\#\{(p_1,p_2,p_3) \, | \, p_i < x, \, p_i \equiv 1 \,(4), \, (\frac{p_j}{p_i})=1, \, [p_1,p_2,p_3]=-1\}}{\pi(x)^3}.
\]
There is an identity
\[
 \begin{split}
 &\#\{(p_1,p_2,p_3) \, | \, p_i < x, \, p_i \equiv 1 \,(4), \, (\frac{p_j}{p_i})=1, \, [p_1,p_2,p_3]=-1\} \\
 &= \sum_{\substack{p_1,p_2<x \\ p_1 \equiv p_2 \equiv 1 \, (4), \, (\frac{p_1}{p_2})=1}}\#\{p \, | \, p < x, \, p \equiv 1 \,(4), \, (\frac{p_1}{p})=(\frac{p_2}{p})=1, \, [p_1,p_2,p]=-1\},
 \end{split}
\]
which is easily derived by the quadratic reciprocity law. Because of this identity, we have
\begin{equation}\label{decomeq}
 \begin{split}
 &\frac{\#\{(p_1,p_2,p_3) \, | \, p_i < x, \, p_i \equiv 1 \,(4), \, (\frac{p_j}{p_i})=1, \, [p_1,p_2,p_3]=-1\}}{\pi(x)^3} \\
 &= \frac{1}{\pi(x)^2}\sum_{\substack{p_1,p_2<x \\ p_1 \equiv p_2 \equiv 1 \, (4), \, (\frac{p_1}{p_2})=1}}\frac{\#\{p \, | \, p < x, \, p \equiv 1 \,(4), \, (\frac{p_1}{p})=(\frac{p_2}{p})=1, \, [p_1,p_2,p]=-1\}}{\pi(x)}.
 \end{split}
\end{equation}

(Step 2): we evaluate the following limit by applying Theorem \ref{EChebo} below:
\begin{equation}\label{plim}
 \lim_{x \to +\infty}\frac{\#\{p \, | \, p < x, \, p \equiv 1 \,(4), \, (\frac{p_1}{p})=(\frac{p_2}{p})=1, \, [p_1,p_2,p]=-1\}}{\pi(x)},
\end{equation}
where $p_1,p_2$ are distinct primes, $p_1 \equiv p_2 \equiv 1 \, (4)$ and $(\frac{p_2}{p_1})=1$. In order to do so, we replace the condition of the numerator of (\ref{plim}) with conditions written by the Artin symbol. \par
Here is the tower of extensions of fields (\ref{extp1p2})
\[
 \begin{tikzcd}
   & k(\sqrt{-1}) \\
 k \arrow[ru, dash]& (k_1k_2)(\sqrt{-1}) \arrow[u, dash]\\
 k_1k_2 \arrow[ru, dash] \arrow[u, dash] & \mathbb{Q}(\sqrt{-1}) \arrow[u, dash]\\
 \mathbb{Q} \arrow[ru, dash] \arrow[u, dash]
\end{tikzcd}
\]
in Section \ref{fieldp1p2}. Note that $k$ is unramified over $\mathbb{Q}$ outside primes $p_1$ and $p_2$. It is easy to see that
\[
 p \equiv 1 \, (4) \hspace{0.2in} \Longleftrightarrow \hspace{0.2in} p \,\, {\rm is} \,\, {\rm completely} \,\, {\rm decomposed} \,\, {\rm in} \,\, \mathbb{Q}(\sqrt{-1}),
\]
\[
 (\frac{p_1}{p})=(\frac{p_2}{p})=1 \hspace{0.2in} \Longleftrightarrow \hspace{0.2in} p \,\, {\rm is} \,\, {\rm completely} \,\, {\rm decomposed} \,\, {\rm in} \,\, k_1k_2.
\]
According to these facts and definition of Borromean primes, it follows that
\[
 \begin{split}
 &[p_1,p_2,p]=-1 \\
 &\Longleftrightarrow \hspace{0.1in} p \,\, {\rm is} \,\, {\rm completely} \,\, {\rm decomposed} \,\, {\rm in} \,\, (k_1k_2)(\sqrt{-1}) \,\, {\rm and} \,\, {\rm not} \,\, {\rm in} \,\, k(\sqrt{-1}).
 \end{split}
\]
For any positive real number $x$, we define sets of primes $S(x), \, S'(x), \, S''(x)$ by
\[
 S(x)=\{p \, | \, p<x, \, p \,\, {\rm is} \,\, {\rm completely} \,\, {\rm decomposed} \,\, {\rm in} \,\, k(\sqrt{-1})\},
\]
\[
 S'(x)=\{p \, | \, p<x, \, p \,\, {\rm is} \,\, {\rm not} \,\, {\rm completely} \,\, {\rm decomposed} \,\, {\rm and} \,\, {\rm unramified} \,\, {\rm in} \,\, k(\sqrt{-1})\}
\]
and
\[
 S''(x)=\{p \, | \, p<x, \, p \,\, {\rm is} \,\, {\rm completely} \,\, {\rm decomposed} \,\, {\rm in} \,\, (k_1k_2)(\sqrt{-1})\}.
\]
Then we have
\begin{equation}\label{4equ1}
 \begin{split}
 \#\{p \, | \, p < x, \, p \equiv 1 \,(4), \, (\frac{p_1}{p})=(\frac{p_2}{p})=1, \, &[p_1,p_2,p]=-1\} \\
  &=\#(S'(x) \cap S''(x)),
 \end{split}
\end{equation}
$S(x) \subset S''(x)$ and if $x>p_1,p_2$, it follows that
\[
 \{p \, | \, p<x\}=S(x) \sqcup S'(x) \sqcup \{p_1,p_2\}
\]
where the operator $\sqcup$ means a disjoint union of sets. Considering an intersection with $S''(x)$ in both sides of the last identity, we see that
\[
 \begin{split}
 S''(x) &= S(x) \sqcup (S'(x) \cap S''(x)).
 \end{split}
\]
Since $S(x)$ and $S'(x) \cap S''(x)$ are disjoint, on the number of elements we obtain
\[
 \#(S'(x) \cap S''(x))=\# S''(x)-\# S(x).
\]
Substituting this to (\ref{4equ1}), the following identity holds:
\[
 \#\{p \, | \, p < x, \, p \equiv 1 \,(4), \, (\frac{p_1}{p})=(\frac{p_2}{p})=1, \, [p_1,p_2,p]=-1\}=\# S''(x)-\# S(x).
\]
By the definition of the Artin symbol, it follows for an unramified prime $p$ and a number field $M$ that
\[
p \,\, {\rm is} \,\, {\rm completely} \,\, {\rm decomposed} \,\, {\rm in} \,\, M \hspace{0.2in} \Longleftrightarrow \hspace{0.2in} \left[\frac{M/\mathbb{Q}}{p}\right]=C({\rm id})=\{{\rm id}\}.
\]
This equivalence leads to the following identities:
\[
 \# S(x) = \#\left\{p \, \middle| \, p<x, \, \left[\frac{k(\sqrt{-1})/\mathbb{Q}}{p}\right]=\{{\rm id}\}\right\}=\pi_{k(\sqrt{-1})/\mathbb{Q}}(x; {\rm id})
\]
and
\[
 \# S''(x) = \#\left\{p \, \middle| \, p<x, \, \left[\frac{(k_1k_2)(\sqrt{-1})/\mathbb{Q}}{p}\right]=\{{\rm id}\}\right\}=\pi_{(k_1k_2)(\sqrt{-1})/\mathbb{Q}}(x; {\rm id}).
\]
Hence we obtain
\begin{equation}\label{limlim}
 \begin{split}
 &\frac{\#\{p \, | \, p < x, \, p \equiv 1 \,(4), \, (\frac{p_1}{p})=(\frac{p_2}{p})=1, \, [p_1,p_2,p]=-1\}}{\pi(x)} \\
 &\hspace{1in}= \frac{\pi_{(k_1k_2)(\sqrt{-1})/\mathbb{Q}}(x; {\rm id})}{\pi(x)}-\frac{\pi_{k(\sqrt{-1})/\mathbb{Q}}(x; {\rm id})}{\pi(x)}.
 \end{split}
\end{equation}\par
Now, we evaluate the limit of RHS of (\ref{limlim}) and the speed of its convergence when $x \to +\infty$ using the following theorem which evaluate the error term of Chebotarev's density theorem for the natural density more precisely under GRH.

\Thm{{\cite[Corollary 1.2.]{GRENIE2019441}}}{EChebo}{Suppose that GRH is true. Let $M$ be a number field, $M'$ be a finite Galois extension of $M$ and $G$ be the Galois group ${\rm Gal}(M'/M)$. For any real number $x \geq 2$,
\[
\begin{split}
&\left|\pi_{M'/M}(x;\sigma)-\frac{\# C(\sigma)}{\# G}\int_2^x\frac{du}{\log u} \right| \\
&\leq \frac{\# C(\sigma)}{\# G} \cdot \sqrt{x}\left[\left(\frac{1}{2\pi}+\frac{3}{\log x}\right)\log \Delta_{M'}+\left(\frac{\log x}{8\pi}+\frac{1}{4\pi}+\frac{6}{\log x}\right)n_{M'} \right]
\end{split}
\]
holds.}

First, we focus on $\frac{\pi_{(k_1k_2)(\sqrt{-1})/\mathbb{Q}}(x; {\rm id})}{\pi(x)}$, the first term of the RHS of (\ref{limlim}). Recall
\[
 [(k_1k_2)(\sqrt{-1}):\mathbb{Q}]=8 \hspace{0.2in} {\rm and} \hspace{0.2in} \Delta_{(k_1k_2)(\sqrt{-1})}=2^8p_1^4p_2^4
\]
(see Proposition \ref{discrilemma} in Section \ref{fieldp1p2}). Applying Theorem \ref{EChebo} with $M=\mathbb{Q}$, $M'=(k_1k_2)(\sqrt{-1})$ and $\sigma={\rm id}$, we have for any primes $p_1,p_2<x$
\[
 \begin{split}
&\left|\pi_{(k_1k_2)(\sqrt{-1})/\mathbb{Q}}(x; {\rm id})-\frac{1}{8}\int_2^x\frac{du}{\log u} \right| \\
&\leq \frac{1}{8}\sqrt{x}\left[\left(\frac{1}{2\pi}+\frac{3}{\log x}\right)\log(2^8p_1^4p_2^4)+\left(\frac{\log x}{8\pi}+\frac{1}{4\pi}+\frac{6}{\log x}\right) \cdot 8 \right] \\
&< \sqrt{x}\left[\left(\frac{1}{2\pi}+\frac{3}{\log x}\right)(\log 2+\log x)+\left(\frac{\log x}{8\pi}+\frac{1}{4\pi}+\frac{6}{\log x}\right)\right].
\end{split}
\]
Dividing the both sides by $\pi(x)$, we see that
\[
 \begin{split}
&\frac{1}{\pi(x)} \cdot \left|\pi_{(k_1k_2)(\sqrt{-1})/\mathbb{Q}}(x; {\rm id})-\frac{1}{8}\int_2^x\frac{du}{\log u} \right| \\
&< \frac{\sqrt{x}}{\pi(x)}\left[\left(\frac{1}{2\pi}+\frac{3}{\log x}\right)(\log 2+\log x)+\left(\frac{\log x}{8\pi}+\frac{1}{4\pi}+\frac{6}{\log x}\right)\right] \\
&= \frac{\frac{x}{\log x}}{\pi(x)} \cdot \frac{(\log x)^2}{\sqrt{x}}\left[\left(\frac{1}{2\pi}+\frac{3}{\log x}\right)\left(\frac{\log 2}{\log x}+1\right)+\left(\frac{1}{8\pi}+\frac{1}{4\pi\log x}+\frac{6}{(\log x)^2}\right)\right].
\end{split}
\]
Because of Prime Number Theorem and ${\displaystyle\lim_{x \to +\infty}}\frac{(\log x)^2}{\sqrt{x}}=0$, there is some positive real number $\alpha$ for any real $\varepsilon>0$ such that for any $x>\alpha$ and any primes $p_1,p_2<x$
\[
 \left|\frac{\pi_{(k_1k_2)(\sqrt{-1})/\mathbb{Q}}(x; {\rm id})}{\pi(x)}-\frac{1}{8 \cdot \pi(x)}\int_2^x\frac{du}{\log u} \right|<\frac{1}{4}\varepsilon.
\]
Therefore some positive real number $\alpha$ exists for any real $\varepsilon >0$ such that for any $x>\alpha$ and any primes $p_1,p_2<x$,
\[
 \begin{split}
 &\left|\frac{\pi_{(k_1k_2)(\sqrt{-1})/\mathbb{Q}}(x; {\rm id})}{\pi(x)}-\frac{1}{8} \right| \\
 &\leq \left|\frac{\pi_{(k_1k_2)(\sqrt{-1})/\mathbb{Q}}(x; {\rm id})}{\pi(x)}-\frac{1}{8 \cdot \pi(x)}\int_2^x\frac{du}{\log u}\right|+\frac{1}{8}\left|\frac{1}{\pi(x)}\int_2^x\frac{du}{\log u}-1 \right| \\
 &<\frac{1}{4}\varepsilon+\frac{1}{8}\left|\frac{1}{\pi(x)}\int_2^x\frac{du}{\log u}-1 \right|,
 \end{split}
\]
where the first inequality is led by the triangle inequality. Because of Prime Number Theorem, there is some positive real number $\alpha'$ such that for any $x>\alpha'$,
\[
 \left|\frac{1}{\pi(x)}\int_2^x\frac{du}{\log u}-1 \right|<2\varepsilon.
\]
Hence we obtain for any $x>\max\{\alpha,\alpha'\}$
\[
 \left|\frac{\pi_{(k_1k_2)(\sqrt{-1})/\mathbb{Q}}(x; {\rm id})}{\pi(x)}-\frac{1}{8} \right| <\frac{1}{4}\varepsilon+\frac{1}{4}\varepsilon=\frac{1}{2}\varepsilon.
\]\par
Furthermore, we also showed
\[
[k(\sqrt{-1}):\mathbb{Q}]=16 \hspace{0.2in} {\rm and} \hspace{0.2in} \Delta_{k(\sqrt{-1})/\mathbb{Q}}=2^{16}p_1^8p_2^8
\]
(see Proposition \ref{discrilemma} in Section \ref{fieldp1p2}). In the manner quite similar to the above argument, we can have the same result for the second term of RHS of (\ref{limlim}). In other words, we can see that there is some real $\alpha>0$ for any real $\varepsilon>0$ such that for any $x>\alpha$ and any primes $p_1,p_2<x$ satisfying $p_1 \equiv p_2 \equiv 1\,(4)$ and $(\frac{p_2}{p_1})=1$
\[
 \left|\frac{\pi_{k(\sqrt{-1})/\mathbb{Q}}(x; {\rm id})}{\pi(x)}-\frac{1}{16} \right| <\frac{1}{2}\varepsilon.
\]\par
Consequently, we can choose some real $\alpha>0$ for any real $\varepsilon >0$ such that for any $x>\alpha$ and any primes $p_1,p_2<x$ satisfying $p_1 \equiv p_2 \equiv 1\,(4)$ and $(\frac{p_2}{p_1})=1$,
\[
 \left|\frac{\pi_{(k_1k_2)(\sqrt{-1})/\mathbb{Q}}(x; {\rm id})}{\pi(x)}-\frac{1}{8} \right| <\frac{1}{2}\varepsilon \hspace{0.2in} {\rm and} \hspace{0.2in} \left|\frac{\pi_{k(\sqrt{-1})/\mathbb{Q}}(x; {\rm id})}{\pi(x)}-\frac{1}{16} \right| <\frac{1}{2}\varepsilon.
\]
Then we have the following because of (\ref{limlim}) and the triangle inequality:
\[
 \begin{split}
 &\left|\frac{\#\{p \, | \, p < x, \, p \equiv 1 \,(4), \, (\frac{p_1}{p})=(\frac{p_2}{p})=1, \, [p_1,p_2,p]=-1\}}{\pi(x)}-\frac{1}{16}\right| \\
 &\leq \left|\frac{\pi_{(k_1k_2)(\sqrt{-1})/\mathbb{Q}}(x; {\rm id})}{\pi(x)}-\frac{1}{8}\right|+\left|\frac{\pi_{k(\sqrt{-1})/\mathbb{Q}}(x; {\rm id})}{\pi(x)}-\frac{1}{16}\right| < \frac{1}{2}\varepsilon+\frac{1}{2}\varepsilon=\varepsilon.
 \end{split}
\]
We show this fact as the next proposition.

\Prp{}{plimthm}{There exists some real $\alpha>0$ for any real $\varepsilon >0$ such that for any $x>\alpha$ and any primes $p_1,p_2<x$ satisfying $p_1 \equiv p_2 \equiv 1\,(4)$ and $(\frac{p_2}{p_1})=1$,
\[
 \left|\frac{\#\{p \, | \, p < x, \, p \equiv 1 \,(4), \, (\frac{p_1}{p})=(\frac{p_2}{p})=1, \, [p_1,p_2,p]=-1\}}{\pi(x)}-\frac{1}{16}\right|<\varepsilon.
\]}

(Step 3): Now, we prove Main Theorem by Proposition \ref{plimthm}. Here is the identity (\ref{decomeq}) again:
\[
 \begin{split}
 &\frac{\#\{(p_1,p_2,p_3) \, | \, p_i < x, \, p_i \equiv 1 \,(4), \, (\frac{p_j}{p_i})=1, \, [p_1,p_2,p_3]=-1\}}{\pi(x)^3} \\
 &= \frac{1}{\pi(x)^2}\sum_{\substack{p_1,p_2<x \\ p_1 \equiv p_2 \equiv 1 \, (4), \, (\frac{p_1}{p_2})=1}}\frac{\#\{p \, | \, p < x, \, p \equiv 1 \,(4), \, (\frac{p_1}{p})=(\frac{p_2}{p})=1, \, [p_1,p_2,p]=-1\}}{\pi(x)}.
 \end{split}
\]
Let $\rho_{p_1,p_2}(x)=\#\{p \, | \, p < x, \, p \equiv 1 \,(4), \, (\frac{p_1}{p})=(\frac{p_2}{p})=1, \, [p_1,p_2,p]=-1\}$ and $(\dagger)$ denote the conditions
\[
 p_1,p_2<x, \, p_1 \equiv p_2 \equiv 1 \, (4), \, (\frac{p_1}{p_2})=1.
\]
Take any real $\varepsilon>0$. According to Proposition \ref{plimthm}, we can choose some real $\alpha>0$ such that for any $x>\alpha$ and any primes $p_1,p_2$ satisfying $(\dagger)$,
\begin{equation}\label{mainineq}
 \left|\frac{\rho_{p_1,p_2}(x)}{\pi(x)}-\frac{1}{16}\right|<\varepsilon.
\end{equation}
According to the inequality (\ref{mainineq}), we obtain for any real number $x>\alpha$
\[
 \begin{split}
 &\left|\frac{1}{\pi(x)^2}\sum_{(p_1,p_2) \in \{(p_1,p_2) \, | (\dagger)\}}\frac{\rho_{p_1,p_2}(x)}{\pi(x)}-\frac{1}{16} \cdot \frac{\#\{(p_1,p_2) \, | \, (\dagger)\}}{\pi(x)^2}\right| \\
 &= \left|\frac{1}{\pi(x)^2}\sum_{(p_1,p_2) \in \{(p_1,p_2) \, | (\dagger)\}}\left(\frac{\rho_{p_1,p_2}(x)}{\pi(x)}-\frac{1}{16}\right)\right| \\
 &\leq \frac{1}{\pi(x)^2}\sum_{(p_1,p_2) \in \{(p_1,p_2) \, | (\dagger)\}}\left|\frac{\rho_{p_1,p_2}(x)}{\pi(x)}-\frac{1}{16}\right| < \frac{\#\{(p_1,p_2) \, | \, (\dagger)\}}{\pi(x)^2} \cdot \varepsilon
 \end{split}
\]
where the second inequality is derived from the triangle inequality. Since $\frac{\#\{(p_1,p_2) \, | \, (\dagger)\}}{\pi(x)^2}$ converge on $\frac{1}{8}<1$ when $x \to +\infty$ because of Theorem \ref{residenthm}, there exists some real $\alpha'>0$ such that for any $x>\alpha'$ we have $\frac{\#\{(p_1,p_2) \, | \, (\dagger)\}}{\pi(x)^2}<1$. Therefore when $x>\max\{\alpha,\alpha' \}$,
\[
 \left|\frac{1}{\pi(x)^2}\sum_{(p_1,p_2) \in \{(p_1,p_2) \, | (\dagger)\}}\frac{\rho_{p_1,p_2}(x)}{\pi(x)}-\frac{1}{16} \cdot \frac{\#\{(p_1,p_2) \, | \, (\dagger)\}}{\pi(x)^2}\right| < \varepsilon
\]
holds. By arbitrariness of $\varepsilon>0$, we see that
\[
\lim_{x \to +\infty}\frac{1}{\pi(x)^2}\sum_{(p_1,p_2) \in \{(p_1,p_2) \, | (\dagger)\}}\frac{\rho_{p_1,p_2}(x)}{\pi(x)}=\lim_{x \to +\infty}\left(\frac{1}{16}\cdot \frac{\#\{(p_1,p_2) \, | \, (\dagger)\}}{\pi(x)^2}\right).
\]
According to Theorem \ref{residenthm}, we obtain
\[
 \lim_{x \to +\infty}\frac{1}{\pi(x)^2}\sum_{(p_1,p_2) \in \{(p_1,p_2) \, | (\dagger)\}}\frac{\rho_{p_1,p_2}(x)}{\pi(x)}=\frac{1}{16} \cdot \frac{1}{8}=\frac{1}{128}.
\]
This completes the proof of Theorem \ref{mainthm}, which is Main Theorem.
\end{proof}

\bibliographystyle{alpha} 
\bibliography{reference}

\vspace{0.5in}

\begin{flushleft}
\footnotesize
Yuki Ishida \\
Graduate School of Mathematics, Kyushu University, \\
744, Motooka, Nishi-ku, Fukuoka, 819-0395, JAPAN \\
e-mail : ishida.yuki.920@s.kyushu-u.ac.jp \\

\vspace{0.1in}

Atsuki Kuramoto \\
Graduate School of Mathematics, Kyushu University, \\
744, Motooka, Nishi-ku, Fukuoka, 819-0395, JAPAN \\
e-mail : kuramoto.atsuki.257@s.kyushu-u.ac.jp\\

\vspace{0.1in}

Dingchuan Zheng \\
Graduate School of Mathematics, Kyushu University, \\
744, Motooka, Nishi-ku, Fukuoka, 819-0395, JAPAN \\
e-mail : zheng.dingchuan.011@s.kyushu-u.ac.jp
\end{flushleft}

\end{document}